\documentclass[12pt,russian,reqno]{amsart}

\usepackage[cp1251]{inputenc}
\usepackage[russian]{babel}

\usepackage{amsmath}
\usepackage{amsfonts}
\usepackage{amssymb}

\theoremstyle{plain}
\newtheorem{theorem}{Теорема}
\newtheorem{hyptheorem}{Гипотеза}

\theoremstyle{remark}
\newtheorem{remark}{Замечание}

\begin{document}

\begin{center}
\textbf{О ДИФФЕРЕНЦИАЛЬНОМ ПРИЗНАКЕ ГОМЕОМОРФНОСТИ ОТОБРАЖЕНИЯ, НАЙДЕННОМ Н.В.~ЕФИМОВЫМ}\footnotemark
\end{center}
\nopagebreak\footnotetext{Работа выполнена при поддержке РФФИ
(код проекта 10--01--91000--анф), ФЦП
<<Научные и научно-педагогические кадры инновационной России>> на 2009--2013 гг. 
(гос. контракт 02.740.11.0457) и Совета по грантам Президента РФ для поддержки ведущих
научных школ (НШ-6613.2010.1).}

\bigskip

\hfill\textbf{В.\,А.~Александров}

\hfill\textit{Институт математики им. С.\,Л.~Соболева СО РАН}

\hfill\textit{и Новосибирский государственный университет}

\hfill\text{e-mail: alex@math.nsc.ru}

\bigskip

\noindent{УДК 514.772, 517.938}

\bigskip

{\small

\begin{center}
\begin{minipage}[t]{140mm}
\textbf{Ключевые слова:} диффеоморфизм, полная  поверхность, гауссова кривизна, гипотеза  якобиана, устойчивость по Ляпунову.
\end{minipage}
\end{center}

\bigskip

\begin{center}
\begin{minipage}[t]{140mm}

\begin{center}
\textbf{Аннотация}
\end{center}

{ }\qquad В 1968 году Н.\,В.~Ефимов доказал следующую замечательную теорему:

{ }\qquad \textit{Пусть $f:\mathbb R^2\to\mathbb R^2\in C^1$, причём
$\det f'(x)<0$ для  всех $x\in\mathbb R^2$.
Пусть, кроме того, существуют положительная  функция  $a=a(x)>0$ и неотрицательные 
постоянные $C_1$, $C_2$ такие, что для  всех $x, y\in\mathbb R^2$ справедливы неравенства
$|1/a(x)-1/a(y)|\leqslant C_1 |x-y|+C_2$ и
$|\det f'(x)|\geqslant a(x)|\mbox{\rm  rot\,}f(x)|+a^2(x).$
Тогда $f(\mathbb R^2)$ есть выпуклая  область и $f$ отображает $\mathbb R^2$ на 
$f(\mathbb R^2)$  гомеоморфно.}

{ }\qquad Здесь $\mbox{\rm rot\,}f(x)$ обозначает ротор функции $f$ в точке 
$x\in\mathbb R^2$.

{ }\qquad Данна  статья   является  обзором аналогов этой теоремы, её обобщений и приложений 
в теории поверхностей, теории функций, а также в исследованиях по гипотезе якобиана и глобальной 
асимптотической устойчивости динамических систем.
\end{minipage}
\end{center}

\begin{center}
\begin{minipage}[t]{140mm}

\begin{center}
\textbf{Abstract}
\end{center}

\textit{Victor Alexandrov, On a differential test of homeomorphism, found by N.\,V.~Efimov.}

{ }\qquad In the year 1968 N.\,V.~Efimov has proven the following remarkable theorem:

{ }\qquad \textit{Let $f:\mathbb R^2\to\mathbb R^2\in C^1$ be such that $\det f'(x)<0$ 
for all $x\in\mathbb R^2$ and let there exist
a function  $a=a(x)>0$ and constants $C_1\geqslant 0$, $C_2\geqslant 0$ 
such that the inequalities
$|1/a(x)-1/a(y)|\leqslant C_1 |x-y|+C_2$
and
$|\det f'(x)|\geqslant a(x)|\mbox{\rm  curl\,}f(x)|+a^2(x)$
hold true for all $x, y\in\mathbb R^2$. 
Then $f(\mathbb R^2)$ is a convex domain and $f$ maps $\mathbb R^2$ onto
$f(\mathbb R^2)$ homeomorhically.}

{ }\qquad Here $\mbox{\rm curl\,}f(x)$ stands for the curl of $f$ at $x\in\mathbb R^2$.

{ }\qquad This article is an overview of analogues of this theorem, its generalizations
and applications in the theory of surfaces, theory of functions, as well as in the study
of the Jacobian conjecture and global asymptotic stability of dynamical systems.
\end{minipage}
\end{center}
}

\section{Введение}

Одним из крупнейших достижений математики XIX века стало создание геометрии Лобачевского 
или, как сейчас принято говорить, --- гиперболической геометрии. 
При этом принципиальное значение имел вопрос о внутренней непротиворечивости такой геометрии.
Первое продвижение в этом направлении сделал в 1868 году Э.~Бельтрами в работе 
<<Опыт интерпретации неевклидовой геометрии>>, 
см. \cite{Be68} или русский перевод в \cite[с.~342--365]{No56}.
Он доказал, что  \textit{специальные\,\footnote{А именно, области одного из следующих трёх типов:
(i) заключённые между двумя пересекающимися прямыми и ортогональной им окружностью, 
(ii) заключённые между двумя расходящимися прямыми, их общим перпендикуляром и 
ортогональной к ним эквидистантой, и 
(iii) заключённые между двумя параллельными прямыми и ортогональным к ним орициклом.}
области плоскости Лобачевского могут быть изометрично 
отображены на подходящие области на псевдосфере,} т.\,е. на поверхности, образуемой 
вращением трактрисы около её асимптоты. Этот результат был воспринят современниками 
Э.~Бельтрами, в сущности, как доказательство непротиворечивости геометрии Лобачевского. 
Впрочем, не в каждой точке псевдосферы существует касательная плоскость и
для  полного успеха нужно было бы найти в трёхмерном евклидовом пространстве регулярную 
поверхность, изометричную всей плоскости Лобачевского.

Однако в 1901 году Д.~Гильберт доказал, что это невозможно, см. \cite{Hi01} 
или русский перевод в \cite[Добавление V]{Hi48}. А именно, он доказал, что 
\textit{в трёхмерном евклидовом пространстве не существует полной регулярной
поверхности с отрицательной и постоянной гауссовой кривизной}. В результате доказательство
непротиворечивости геометрии Лобачевского пошло другим путём, однако вопрос о том насколько
существенно требование постоянства гауссовой кривизны в теореме Д.~Гильберта долгое время  
оставался  открытым.

Только в 1964 году в Н.\,В.~Ефимов полностью снял это ограничение, доказав,
что \textit{в трёхмерном евклидовом пространстве невозможна полна  регулярна  поверхность с 
гауссовой кри­визной $K\leqslant\mbox{\rm const\,}<0$}, см. \cite{Ef64} или \cite{KM86}.
Для  доказательства он разработал оригинальный метод исследования сферического отображения,
называемый иногда <<методом вогнутой опоры>>, с помощью которого
в 1968 году доказал следующие две теоремы:

\begin{theorem}[Н.\,В.~Ефимов \cite{Ef68}]\label{Th_II_Ef68}
Пусть $f:\mathbb R^2\to\mathbb R^2$ принадлежит классу $C^1$, 
причём якобиан отображения  $f$ всюду отрицателен, т.\,е.
$\det f'(x)<0$ для  всех $x\in\mathbb R^2$.
Пусть, кроме того, существуют положительная  функция  $a=a(x)>0$ и неотрицательные 
постоянные $C_1$, $C_2$ такие, что для  всех $x, y\in\mathbb R^2$ справедливо неравенство 
\begin{equation}\label{Eq_1}
\biggl|\dfrac{1}{a(x)}-\dfrac{1}{a(y)}\biggr|\leqslant C_1 |x-y|+C_2.              
\end{equation}
Тогда если для  всех $x\in\mathbb R^2$ выполнено неравенство
$$
|\det f'(x)|\geqslant a(x)|\mbox{\rm  rot\,}f(x)|+a^2(x),           
$$
то $f(\mathbb R^2)$ есть выпуклая  область и $f$ отображает $\mathbb R^2$ на 
$f(\mathbb R^2)$  гомеоморфно.
\end{theorem}

\begin{remark}
Здесь $\mbox{\rm rot\,}f(x)$ означает ротор функции $f$ в точке 
$x=(x_1, x_2)\in\mathbb R^2$, т.\,е. 
$$\mbox{\rm rot\,}f(x)= \frac{\partial f_2}{\partial x_1}(x)-\frac{\partial f_1}{\partial x_2}(x).$$
\end{remark}

\begin{remark}
Если функция  $1/a(x)$ имеет экспоненциальный рост (так что усло­вие (\ref{Eq_1}) 
грубо не соблюдается), то утверждение теоремы \ref{Th_II_Ef68} даже для  потен­циальных отображений
(т.\,е. таких, что $\mbox{\rm rot\,}f(x)\equiv 0$) может не иметь места. 
Например, для  отображения 
$f(x_1,x_2)=(e^{x_1}\sin x_2, e^{x_1}\cos x_2)$ 
имеем $\mbox{\rm rot\,}f(x)\equiv 0$, $\det f'(x_1,x_2) = - e^{2x_1}$, 
$f(\mathbb R^2)=\mathbb R^2\diagdown \{(0,0)\}$, и каждая  точка из $f(\mathbb R^2)$
имеет бесконечное число прообразов.
\end{remark}

\begin{theorem}[Н.\,В.~Ефимов \cite{Ef68}]\label{Th_3_Ef68}
Пусть $f:\mathbb R^2\to\mathbb R^2$ принадлежит классу $C^1$
и пусть существует постоянна  $a= \mbox{\rm const\,} > 0$, 
такая, что для всех $x\in\mathbb R^2$ выполнено неравенство
$$
|\det f'(x)|\geqslant a|\mbox{\rm  rot\,}f(x)|+a^2.
$$
Тогда $f(\mathbb R^2)$ есть либо бесконечная  полоса между параллельными прямыми,  либо полуплоскость,
либо плоскость и $f$ отображает $\mathbb R^2$ на $f(\mathbb R^2)$ гомеоморфно.
\end{theorem}

\begin{remark}
Условия  теоремы \ref{Th_3_Ef68} заведомо выполнены, если  для  всех $x\in\mathbb R^2$
справедливы неравенства
$|\det f'(x)|\geqslant \mbox{\rm const}> 0$ 
и $|\mbox{\rm  rot\,}f(x)|\leqslant \mbox{\rm const}$.
\end{remark}

В настоящей статье речь идёт о следе, оставленном теоремами \ref{Th_II_Ef68} и \ref{Th_3_Ef68} 
в математике. Точнее, мы даём обзор аналогов этих теорем, их обобщений и приложений, 
за последние 40 лет.
Отдельные разделы посвящены обзору работ, мотивированных исследованиями по теории поверхностей, 
теории функций, 
исследованиям по гипотезе якобиана и 
глобальной асимптотической устойчивости динамических систем.

\section{Результаты, мотивированные теорией поверхностей}\label{TS}

В статьях Б.\,Е.~Кантора \cite{Ka70} и С.\,П.~Гейсберга \cite{Ge70}, опубликованных в 1970 году
и посвящённых изучению нормального образа полной поверхности отрицательной кривизны, 
исследован, в частности, вопрос о том, какая из указанных Н.\,В.~Ефимовым в 
теореме~\ref{Th_3_Ef68} возможностей для $f(\mathbb R^2)$ в действительности реализуется. 
Опишем эти результаты, вошедшие в учебник \cite[\S\,31]{BVK73}, подробнее.

Для линейного отображения $f_1(x_1,x_2)=x_2$,
$f_2(x_1,x_2)=x_1$ область $f(\mathbb R^2)$  является плоскостью.
Так что в условиях теоремы~\ref{Th_3_Ef68} область $f(\mathbb R^2)$ действительно может быть плоскостью. 

В \cite{Ka70} Б.\,Е.~Кантор привёл пример следующего отображения $f:\mathbb R^2\to\mathbb R^2$
\begin{equation}\label{Eq_Ka70}
f_1(x_1,x_2) = \mbox{ln\,}\bigl(x_1+\sqrt{x_1^2+ e^{-2x_2}}\bigr) + x_2,\qquad
f_2(x_1,x_2) = \sqrt{x_1^2 + e^{-2x_2}},
\end{equation}
для которого 
\begin{equation}\label{Eq_rot=0}
\det f'(x)\equiv \mbox{\rm const}<0 \quad\mbox{и}\quad \mbox{\rm rot\,}f(x)\equiv 0,
\end{equation} 
и, вместе с тем, $f(\mathbb R^2)$ есть полуплоскость. Так что и эта возможность, указанная в 
теореме~\ref{Th_3_Ef68}, реализуется. 

Основная идея Б.\,Е.~Кантора состоит в том, что условия (\ref{Eq_rot=0})
будут заведомо выполнены, если отображение $f=(f_1,f_2)$ потенциально (т.\,е. если
существует функция $\varphi:\mathbb R^2\to\mathbb R$, называема  потенциалом, такая, что 
$f_j=\partial\varphi/\partial x_j$ для  $j=1,2$)  и потенциал удовлетворяет уравнению Монжа -- Ампера
$$
\frac{\partial^2 \varphi}{\partial x_1^2}\frac{\partial^2 \varphi}{\partial x_2^2}-
\biggl(\frac{\partial^2 \varphi}{\partial x_1\partial x_2}\biggr)^2=-1.
$$
Далее он использует классическое параметрическое представление для
решений уравнения Монжа -- Ампера \cite{Go36} и получает не одно отображение
(\ref{Eq_Ka70}), а целый класс нетривиальных примеров отображений $f:\mathbb R^2\to\mathbb R^2$,
удовлетворяющих условиям (\ref{Eq_rot=0}), для которых $f(\mathbb R^2)$ есть полуплоскость. 

В \cite{Ka70} Б.\,Е.~Кантор доказал также, что \textit{если
выполнены условия} (\ref{Eq_rot=0}), 
\textit{то $f(\mathbb R^2)$ не может быть полосой}
и высказал гипотезу о том, что 
\textit{в условиях теоремы  \ref{Th_3_Ef68} Н.В.~Ефимова 
область $f(\mathbb R^2)$ вообще не может быть полосой.}
Эта гипотеза остаётся недоказанной и поныне.

С.\,Л.~Гейсберг в \cite{Ge70} несколько усилил результат Б.Е.~Кантора,
показав, что \textit{$f(\mathbb R^2)$ не может быть полосой ни в одном из следующих
случаев}:                                         

(а) \textit{всюду вне некоторого круга выполнены условия 
$\mbox{\rm rot\,}f(x)\equiv 0$, $\det f'(x)=-g^2(x)$, 
$g(x)\geqslant\mbox{\rm const}>0$, причём функция $g$ выпукла вниз};

(б) \textit{$\mbox{\rm rot\,}f(x)\equiv 0$, а $\det f'(x)$ есть многочлен, 
принимающий только отрицательные значения, равномерно отделённые от нуля.}

Метод, применяемый С.\,Л.~Гейсбергом, близок оригинальному методу Н.В.~Ефимова и
основан на изучении характеристик отображения 
$f:\mathbb R^2\to\mathbb R^2$, т.\,е. линий на плоскости $(x_1, x_2)$, на которых
выполняются либо соотношения 
$df_1 = g(x_1,x_2)dx_2$ и $df_2 = -g(x_1,x_2)dx_1$,
либо соотношения
$df_1 = -g(x_1,x_2)dx_2$ и $df_2 = g(x_1,x_2)dx_1$.

Поскольку отображение $g$, обратное к потенциальному отображению $f$, само является потенциальным,
то теорему~\ref{Th_3_Ef68} можно переформулировать следующим образом:
\textit{Пусть область $D\subset\mathbb R^2$ такова, что существует
потенциальный обращающий ориентацию $C^1$-диффеоморфизм $g$, отображающий область
$D$ на всю плоскость $\mathbb R^2$ и
такой, что для всех $y\in\mathbb R^2$  и некоторой постоянной $b= \mbox{\rm const\,} > 0$
выполнено неравенство $|\det g'(b)|\leqslant b$. 
Тогда $D$ есть либо бесконечная  полоса между параллельными прямыми,  либо полуплоскость,
либо плоскость.}

В диссертации Г.\,Я.~Перельмана \cite{Pe90}, посвящённой исследованию седловых поверхностей 
в евклидовых пространствах, изучен, в частности, вопрос о том, какими могут быть области $D$,
если в предыдущем утверждении отказаться от требования ограниченности якобиана диффеоморфизма $g$. 
А именно, им получена следующая теорема:

\begin{theorem}[Г.\,Я.~Перельман {\cite[с.~78]{Pe90}}]\label{Th_Pe90}
Выпуклые области, допускающие потенциальный обращающий ориентацию $C^1$-диффеоморфизм 
на плоскость --- это в точности области из следующего списка: 

{\rm (1)} треугольники{\rm ;}

{\rm (2)} четырёхугольники{\rm ;}

{\rm (3)} плоскость{\rm ;}

{\rm (4)} полуплоскость{\rm ;}

{\rm (5)} полоса{\rm ;}

{\rm (6)} угол{\rm ;}

{\rm (7)} полуполоса{\rm ;}

{\rm (8)} угол с отрезанной вершиной{\rm ;}

{\rm (9)} полуполоса с отрезанной вершиной.
\end{theorem}

Доказательство теоремы~\ref{Th_Pe90}, предложенное Г.\,Я.~Перельманом, является модификацией
доказательства теоремы~\ref{Th_3_Ef68}, данного Н.\,В.~Ефимовым.

\section{Результаты, мотивированные теорией функций}\label{TF}

В работах \cite{Al90, Al91b} показано, что следующая теорема 
эквивалентна теореме~\ref{Th_II_Ef68} Н.\,В.~Ефимова.

\begin{theorem}[В.\,А.~Александров \cite{Al90, Al91b}]\label{Th_Al90}
Пусть $f:\mathbb R^2\to\mathbb R^2$ принадлежит классу $C^1$, 
причём $\det f'(x)<0$ для  всех $x\in\mathbb R^2$. 
Пусть, кроме того, существуют положительная  функция  $a(x)>0$ и неотрицательные 
постоянные $C_1, C_2$ такие, что для  всех $x,y\in\mathbb R^2$ справедливо неравенство
$$
\biggl|\frac{1}{a(x)}-\frac{1}{a(y)}\biggr|\leqslant C_1 |x-y| +C_2.
$$
Тогда если для  всех $x\in\mathbb R^2$ выполнено неравенство
$$
|\mu_2(x)|\geqslant |\mu_1(x)|\geqslant a(x),       
$$
где $\mu_1(x)$ и $\mu_2(x)$ --- собственные числа линейного отображения  $f'(x)$, то
$f(\mathbb R^2)$ есть выпуклая  область и $f$ отображает $\mathbb R^2$ на $f(\mathbb R^2)$ 
гомеоморфно.
\end{theorem}

Предложенное в \cite{Al90, Al91b} доказательство теоремы \ref{Th_Al90} состоит в том,
чтобы убедиться, что из условий теоремы \ref{Th_Al90} вытекают условия теоремы \ref{Th_II_Ef68}
Н.\,В.~Ефимова.
Как указано в~\cite{Al90, Al91b}, Н.\,С.Даирбековым было замечено что и наоборот 
из условий теоремы~\ref{Th_II_Ef68} вытекают условия теоремы~\ref{Th_Al90}. 
В этом смысле теоремы~\ref{Th_II_Ef68} и~\ref{Th_Al90} эквивалентны.

Теорема~\ref{Th_Al90}, по нашему мнению, указывает направление, в котором следует 
искать многомерные аналоги теоремы~\ref{Th_II_Ef68} Н.\,В.~Ефимова. 
Для  б\'ольшей  ясности сформулируем простейшую гипотезу, лежащую в этом направлении.

\begin{hyptheorem}[В.\,А.~Александров \cite{Al90, Al91b}]\label{Th_hyp_Al90}
Пусть $n\geqslant 2$ и $f:\mathbb R^n\to\mathbb R^n$ 
принадлежит классу $C^1$, причём $\det f'(x)\neq 0$ для всех $x\in\mathbb R^n$. 
Пусть, кроме того, существует постоянная $K<+\infty$ такая, что для всех 
$x\in\mathbb R^n$ справедливо неравенство 
$$
\rho(f'(x)^{-1})\leqslant K,
$$
где $\rho(f'(x)^{-1})$ --- спектральный радиус линейного отображения  $f'(x)^{-1}$,
т.\,е. максимум модуля  его собственных чисел. 
Тогда $f(\mathbb R^n)$  является выпуклой областью и отображение $f$ инъективно.
\end{hyptheorem}

Заметим, что, как следует из теоремы~\ref{Th_Al90}, 
при $n=2$ и $\det f'(x)<0$ гипотеза~\ref{Th_hyp_Al90} верна.

В теории функций хорошо известен цикл результатов, обычно называемых
<<теоремой о глобальной обратной функции>>. 
По-видимому, первый результат такого рода был получен Ж.~Адамаром~\cite{Ha06}
в 1906 году для отображений $\mathbb R^2$ в себя. 
Со временем выяснилось, что подобного рода утверждения справедливы и
для отображений банаховых пространств.
Для определённости мы сформулируем один из этих результатов:

\begin{theorem}[Р.~Пласток \cite{Pl74}]\label{Th_Pl74}
Пусть $b$ и $B$ --- банаховы пространства, а непрерывно дифференцируемое отображение
$f:b\to B$ является локальным гомеоморфизмом, причём линейное отображение $f'(x)$ 
обратимо при любом $x\in b$ и 
\begin{equation}\label{Eq_4}
\int\limits_{0}^{+\infty} \inf\limits_{
\begin{array}{c}
\scriptstyle \| x\|\leqslant t \\ \scriptstyle x\in b
\end{array}} 
\| f'(x)^{-1}\|^{-1}\, dt=+\infty.
\end{equation}
Тогда $f$ взаимно однозначно отображает $b$ на всё пространство $B$.
\end{theorem}

Отметим, что если  существует постоянная $K<+\infty$ такая, что для всех 
$x\in b$ справедливо неравенство 
\begin{equation}\label{Eq_5}
\| f'(x)^{-1}\|\leqslant K,
\end{equation}
то условие~(\ref{Eq_4}) теоремы~\ref{Th_Pl74} заведомо выполнено.
А поскольку $\rho(f'(x)^{-1})\leqslant \| f'(x)^{-1}\|$, то
для отображений пространства $\mathbb R^n$, удовлетворяющих условию~(\ref{Eq_5}),
гипотеза~\ref{Th_hyp_Al90} заведомо верна. Более того, для таких отображений
$f(\mathbb R^n)$ не просто является выпуклой областью, но совпадает с $\mathbb R^n$.

Укажем ещё один естественный вопрос, связанный с гипотезой~\ref{Th_hyp_Al90},
остающийся открытым до сих пор. Он состоит в том, чтобы
\textit{найти оценку радиуса шара, где существует обратное отображение $f^{-1}$, 
в терминах спектрального радиуса $\rho(f'(x)^{-1})$,} аналогичную
оценке из следующей теоремы, обобщающей теорему~\ref{Th_Pl74}:

\begin{theorem}[Ф.~Джон \cite{Jo68}]\label{Th_Jo68}
Пусть $b$ и $B$ --- банаховы пространства, 
$\omega\subset b$ --- шар радиуса $r$ с центром $x$, 
$f: \omega \to B$ --- локальный $C^1$-гомеоморфизм.
Тогда обратное отображение $f^{-1}$ существует в шаре $\Omega\subset B$ с центром $f(x)$ и радиусом
$$
R=\int\limits_{0}^{r} \inf\limits_{
\begin{array}{c}
\scriptstyle \| x-y\|\leqslant t \\ \scriptstyle y\in \omega
\end{array}} 
\| f'(y)^{-1}\|^{-1}\, dt.
$$
\end{theorem}

\section{Результаты, мотивированные гипотезой якобиана}\label{JC}

Знаменитая гипотеза якобиана утверждает, что
\textit{если $f: \mathbb C^n\to\mathbb C^n$ полиномиально и 
$\det f'(x)\neq 0$ для всякого $x\in\mathbb C^n$, то $f$ имеет полиномиальное обратное.}
Здесь $x=(x_1,\dots,x_n)\in \mathbb C^n$,
$f(x)=(f_1(x),\dots,f_n(x))$ и каждая компонента
$f_j(x)$ является многочленом от $n$ переменных $x_1,\dots,x_n$.

Гипотезе якобиана посвящена обширная литература, см. например, книгу~\cite{Es00}
и указанные в ней ссылки. 
Однако в приведёном выше виде она не доказана до сих пор.
О её значении говорит, хотя бы, тот факт, что С.~Смейл включил её в свой
знаменитый список <<Математические проблемы для следующего века>>~\cite{Sm00}.

Отметим, что известны различные дополнительные условия на отображение $f$,
при выполнении которых гипотеза якобиана заведомо справедлива.
Одним из таких дополнительных условий является требование инъективности
отображения $f$.\footnote{Элементарное доказательство того, что если $f$ инъективно, то
оно сюръективно и $f^{-1}$ полиномиально, дано в~\cite{Ru95}.}
Этим объясняется интерес специалистов по проблеме якобиана к 
условиям, гарантирующим инъективность отображений. Укажем одну из 
высказанных ими гипотез:

\begin{hyptheorem}[М.~Чемберлэнд~\cite{CM98}]\label{Th_hyp_CM98}
Пусть $n\geqslant 1$ и $f:\mathbb R^n\to\mathbb R^n$ 
принадлежит классу $C^1$, причём $\det f'(x)\neq 0$ для всех $x\in\mathbb R^n$. 
Пусть, кроме того, существует постоянная $K<+\infty$ такая, что для всех 
$x\in\mathbb R^n$ справедливо неравенство 
$$
\rho(f'(x)^{-1})\leqslant K,
$$
где $\rho(f'(x)^{-1})$ --- спектральный радиус линейного отображения  $f'(x)^{-1}$.
Тогда отображение $f$ инъективно.
\end{hyptheorem}

Подчеркнём, что гипотеза~\ref{Th_hyp_CM98} в целом совпадает 
с гипотезой~\ref{Th_hyp_Al90}, хотя формально и слабее последней. 
Но и гипотеза~\ref{Th_hyp_CM98} на сегодня остаётся открытой. 
Наиболее близким к ней по духу результатом является следующая
теорема:

\begin{theorem}[М.~Чемберлэнд, Г.~Мейстерс~\cite{CM98}]\label{Th_CM98} 
Пусть $n\geqslant 1$ и непрерывно дифференцируемое отображение
$f:\mathbb R^n\to\mathbb R^n$ таково, что $\det f'(x)\neq 0$ для всех $x\in\mathbb R^n$.
Пусть, кроме того, для всех $x\in\mathbb R^n$ выполняется неравенство
$$\rho\biggl(\bigl(f'(x)f'(x)^*\bigr)^{-1}\biggr)\leqslant \, \text{\rm const\,},$$
где $\rho (L)$~--- спектральный радиус линейного отображения $L$, а  
$L^*$~--- отображение, сопряжённое к $L$.
Тогда $f$ инъективно.
\end{theorem}

\section{Результаты, мотивированные изучением глобальной асимптотической устойчивости динамических систем}\label{DS}

Проблема глобальной асимптотической устойчивости динамических систем на плоскости
может быть сформулирована так:
\textit{Рассмотрим автономную систему на плоскости
\begin{equation}\label{Eq_6}
\begin{cases}
\dot x_1 &= f_1(x_1,x_2),\\
\dot x_2 & = f_2(x_1,x_2),
\end{cases}
\end{equation}
для которой точка $(x_1,x_2)=(0,0)$ является особой\footnote{Т.\,е. $f_1(0,0)=f_2(0,0)=0$.}, 
а собственные числа матрицы Якоби отображения $f=(f_1,f_2)$ 
имеют отрицательные вещественные части в любой точке плоскости.}

\textit{Верно ли, что тогда тривиальное решение $(x_1,x_2)=(0,0)$ 
глобально асимптотически устойчиво, т.\,е. верно ли, что любое решение
стремится к точке $(0,0)$ при $t\to\infty$?}

Проблема глобальной асимптотической устойчивости играет важную роль в теории и приложениях
динамических систем, см., например, \cite{Ch06}. 

Согласно знаменитой теореме Ляпунова, автономная система линейных дифференциальных уравнений
$\dot{\boldsymbol{x}}=A\boldsymbol{x}$, $ \boldsymbol{x}\in\mathbb R^n$,
является асимптотически устойчивой, если все собственные значения постоянной матрицы~$A$
имеют отрицательные действительные части.
Таким образом, проблема глобальной асимптотической устойчивости отвечает на вопрос о том
можно ли делать вывод об асимптотической устойчивости любого решения нелинейной системы~(\ref{Eq_6}),
как только установлена асимптотическая устойчивость любого решения каждой линеаризованной системы.

В 1963 году Ч.~Олех обнаружил тесную связь между проблемой глобальной асимптотической устойчивости 
и инъективностью некоторых отображений. А именно, он доказал следующую теорему:

\begin{theorem}[Ч.~Олех~\cite{Ol63}]\label{Th_Ol63}  
Следующие два утверждения эквивалентны:

{\rm (1)} проблема глобальной асимптотической устойчивости имеет положительное решение{\rm ;}

{\rm (2)} всякое непрерывно дифференцируемое отображение $f:\mathbb R^2\to\mathbb R^2$, 
собственные числа матрицы Якоби которого всюду имеют отрицательные вещественные части, 
с необходимостью инъективно.
\end{theorem}

В 1994--1995 годах положительное решение проблемы глобальной асимптотической устойчивости
почти одновременно, но независимо друг от друга, получили 
А.\,А.~Глуцюк~\cite{Gl94}, Р.~Фесслер~\cite{Fe95} и К.~Гутьеррес~\cite{Gu95}.
Все три решения опираются на теорему Олеха. 
При этом интерес к проблеме глобальной асимптотической устойчивости не исчезает,
напротив, --- появляются всё новые её доказательства, см., например,~\cite{CHQ01}.
Более того, специалисты по глобальной асимптотической устойчивости
получили разнообразные дифференциальные признаки инъективности отображений,
родственные теоремам~\ref{Th_II_Ef68} и ~\ref{Th_3_Ef68}, и
лежащие в русле гипотез~\ref{Th_hyp_Al90} и \ref{Th_hyp_CM98}.
Приведём здесь лишь два из них.  Для этого обозначим через
$\text{Spec\,} (f)$ совокупность всех собственных значений 
линейных отображений $f'(x)$ для всех $x\in\mathbb R^2$.

\begin{theorem}[К.~Гутьеррес~\cite{Gu95}]\label{Th_Gu95}  
Пусть непрерывно дифференцируемое отображение $f:\mathbb R^2 \to\mathbb R^2$ таково, что
$\text{\rm Spec\,} (f)\cap [0,+\infty)= \varnothing.$
Тогда $f$ инъективно.
\end{theorem}

\begin{theorem}[М.~Кобо, К.~Гутьеррес, Ж.~Ллибре~\cite{CGL02}]\label{Th_CGL02}  
Пусть $f:\mathbb R^2 \to\mathbb R^2$ --- непрерывно дифференцируемое отображение и пусть
существует $\varepsilon >0$ такое, что
$\text{\rm Spec\,} (f)\cap(-\varepsilon, \varepsilon)= \varnothing.$
Тогда $f$ инъективно. 
\end{theorem}

Связь этих теорем с гипотезами~\ref{Th_hyp_Al90} и \ref{Th_hyp_CM98},
а значит, и с теоремами~\ref{Th_II_Ef68} и ~\ref{Th_3_Ef68} Н.\,В.~Ефимова, очевидна.
Например, используя обозначение $\text{Spec\,} (f)$ гипотезу~\ref{Th_hyp_CM98}
можно сформулировать так: 
\textit{Если непрерывно дифференцируемое отображение $f:\mathbb R^n\to\mathbb R^n$
таково, что существует  $\varepsilon >0$ такое, что круг с центром в нуле и радиусом  
$\varepsilon >0$ на комплексной плоскости не пересекается с множеством $\text{\rm Spec\,} (f)$,
то $f$ инъективно.}
Таким образом, при $n=2$, теорема~\ref{Th_CGL02} показывает, что справедливо даже
более сильное, чем гипотеза~\ref{Th_hyp_CM98}, утверждение: в этом случае для инъективности $f$ 
достаточно потребовать, чтобы $\text{\rm Spec\,} (f)$ не пересекался только с достаточно малым
интервалом $(-\varepsilon, \varepsilon)$ вещественной оси, а не с целым кругом 
радиуса $\varepsilon$ с центром в нуле.

\section{Заключение}

Данная статья является записью доклада, прочитанного автором на Международной конференции 
<<Метрическая геометрия поверхностей и многогранников>>, посвященной 100-летию со дня рождения 
Николая Владимировича Ефимова и его научному наследию, проходившей с 18-го по 21-е августа 2010 года
на механико-математическом факультете Московского государственного университета имени 
М.\,В.~Ломоносова.

В статье упомянуты все известные нам области математики, связанные с теоремами~\ref{Th_II_Ef68}
и \ref{Th_3_Ef68} Н.\,В.~Ефимова. Но, будучи ограничены во времени, зачастую мы не могли
в пределах каждой области перечислить все результаты, связанные с теоремами~\ref{Th_II_Ef68}
и \ref{Th_3_Ef68}. В такой ситуации мы приводили только самые сильные в своей 
области результаты, доступные нам на момент написания статьи. 
Следует иметь ввиду, что зачастую эти результаты были получены неоднократным улучшением
более ранних и менее общих теорем. Ссылки на последние заинтересованный читатель может 
найти в цитируемых нами работах.

Однако даже из нашего краткого изложения видно, что в теоремах~\ref{Th_II_Ef68}
и \ref{Th_3_Ef68} Н.\,В.~Ефимов открыл глубокий математических факт,
который связан с разными разделами математики, который непосредственно стимулировал исследования
одних математиков, либо был переоткрыт заново другими, не знакомыми с работами Н.\,В.~Ефимова.
Этот факт прочно связан с именем Николая Владимировича Ефимова и продолжает 
до сих пор интересовать математиков.

\bigskip
\small{
{}\hfill\textit{Статья поступила 18 октября 2010 г.}}

\end{document}